\documentclass[10pt, a4paper]{amsart}
\usepackage{amsfonts}
\usepackage{amsthm}
\usepackage{amsmath}

\theoremstyle{plain} \newtheorem{lemma}{Lemma}[section]
\theoremstyle{plain} \newtheorem{proposition}[lemma]{Proposition}
\theoremstyle{plain} \newtheorem{theorem}[lemma]{Theorem}
\theoremstyle{plain} \newtheorem{corollary}[lemma]{Corollary}

\theoremstyle{definition}
\newtheorem{algorithm}{Algorithm}
\newtheorem{remark}[lemma]{Remark}

\usepackage[all]{xy}

\newcommand{\Q}{\mathbb{Q}}
\newcommand{\C}{\mathbb{C}}
\newcommand{\Z}{\mathbb{Z}}
\newcommand{\GL}{\mathop{\mathrm{GL}}}
\newcommand{\g}{\mathop{\mathfrak{g}}}
\newcommand{\h}{\mathop{\mathfrak{h}}}
\newcommand{\gl}{\mathop{\mathfrak{gl}}}

\begin{document}

\title{Constructing arithmetic subgroups of unipotent groups}
\author{Willem A. de Graaf}
\author{Andrea Pavan}
\thanks{The second author would like to thank the University of Trento
for its hospitality during 2008, when the work on the paper was carried out.}
\address{Dipartimento di Matematica\\
Universitry of Trento\\
Via Sommarive 14\\
38050 Povo (Trento)\\
Italy}
\email{degraaf@science.unitn.it}
\address{Dipartimento di Matematica Pura ed Applicata\\
University of Padova\\
Via Trieste 63 \\
35121 Padova\\
Italy}
\email{pan@math.unipd.it}

\date{}
\maketitle

\begin{abstract}
Let $G$ be a unipotent algebraic subgroup of some $\textrm{GL}_{m}(\mathbb{C})$
defined over $\mathbb{Q}$. We describe an algorithm for finding a finite set 
of generators of the subgroup $G(\Z) = G \cap \textrm{GL}_{m}
(\mathbb{Z})$. This is based on a new proof of the result (in more general
form due to Borel and Harish-Chandra) that such a finite generating set 
exists. 
\end{abstract}

\section{Introduction}

Let $G$ be an algebraic subgroup of $\textrm{GL}_{m}(\mathbb{C})$ defined over 
$\mathbb{Q}$, where $m \geq 1$. Then for a subring $R$ of $\C$ we set
\begin{displaymath}
G(R) = G \cap \textrm{GL}_{m}(R).
\end{displaymath}
The group $G(\Z)$ and any other subgroup $\Gamma$ of $G(\Q)$ commensurable 
with it are called arithmetic subgroups of $G$. 

Arithmetic groups occur in many contexts. Examples are: the automorphism group
of a finitely generated nilpotent group (\cite{segal}, Chapter 6), the group 
of units of the ring of integers of a number field, and the group of units
of the group algebra $\Z \mathcal{G}$, where $\mathcal{G}$ is a finite group. 
A celebrated theorem 
of Borel and Harish-Chandra (\cite{borelhs}) says that arithmetic groups are 
finitely generated.
In this paper we consider the problem of computing a finite set of generators
of an arithmetic subgroup of a unipotent group. \par
This problem was also treated in the paper \cite{grunewald_segal} by 
Grunewald and Segal, where a general 
algorithm for all arithmetic groups was outlined. However, their declared aim 
was to show that such a computation is, at least in principle, feasible, and 
no attempt was made to make the algorithms as efficient as possible. 
And unfortunately it 
appears to be extremely hard to use their algorithm in practice.

In this paper we describe a practical algorithm for finding a finite set
of generators of $G(\Z)$, in case 
$G$ is unipotent, that is to say, all its elements are unipotent matrices. 
As a byproduct this yields an independent 
proof of the Borel-Harish-Chandra theorem in this case. Also, we can show that 
the groups $G_L$ are $T$-groups of Hirsh length equal to $\dim G$. In order 
to show that the algorithm is practical we have implemented it in the language 
of the  computer algebra system {\sf GAP}4 (cf. \cite{gap4}). 

We now sketch the main idea of the algorithm. Let $V$ be the vector space
on which $G$ acts naturally. Let 
$$0 = V_{0} < V_{1} < \cdots < V_{n} = V$$
be a flag of $V$ with respect to the action of $G$ (this means that 
for $v\in V_i$ and $g\in G$ we have $gv \equiv v \bmod V_{i-1}$). Then 
we can form the $G$-module $V^\star = V_{n-1}\oplus \frac{V}{V_1}$. In informal
terms the matrix of a $g\in G$ acting on $V^\star$ is formed from 
the matrix of its action on $V$ by taking the block in the upper left part
of the matrix, and the block in the bottom right part of the matrix, and 
constructing the block matrix consisting of these two blocks.
Now let $Q$ be the image of $G$ in $\GL(V^\star)$; then we can recursively
compute generators of $Q(\Z)$. The recursion works because the $Q$-flag
in $V^\star$ has smaller length. Let $\pi : G\to Q$ be the 
projection. In Section \ref{thresults} we describe $\pi(G(\Z))$ (Proposition
\ref{kerpsi}). We show how to find generators of $\pi(G(\Z))$, and their preimages
in $G(\Z)$. Let $N(\Z)\subset G(\Z)$ denote the kernel of $\pi$. We find a finite
set of generators of $N(\Z)$, and joined to the elements of $G(\Z)$ found 
earlier this solves the problem. \par
These ideas are detailed in Section \ref{thresults}. In Section 
\ref{sect:exa} we illustrate them with a simple example. 
The constructions of Section \ref{thresults} do not immediately
yield an implementable algorithm. In order to obtain that we need some 
technical preparation.
In Section \ref{sect:basis} we describe some results that allow us to work
with the Lie algebra of $G$ rather than with $G$ itself. Section 
\ref{sect:tgrps} contains some material on $T$-groups. In Section 
\ref{sect:smith} we describe some algorithms for lattices that
we need. Then in Section \ref{sect:main} we give a detailed description of the
main algorithm, and prove its correctness. The last section
describes some practical experiences with our implementation of this
algorithm in {\sf GAP}4.

\section{The derived representation} \label{thresults}

The goal of this section is to introduce some notation, and to prove
the results that underpin the main algorithm.\par
Let $V$ be a finite dimensional vector space over $\mathbb{Q}$, and $L$ a 
full-dimensional lattice of $V$. We do not prove the next lemma here; it
will follow from Lemma \ref{lem:Lcomp}.
\begin{lemma}
Let $U$ and $W$ be two subspaces of $V$  with $U \subseteq W$. Then there exist subspaces $U'$ and $W'$ of $V$ such that 
\begin{displaymath}
W' \subseteq U'
\end{displaymath}
and equalities
\begin{displaymath}
U \oplus U' = V = W \oplus W'
\end{displaymath}
and
\begin{displaymath}
\big ( L \cap U \big ) + \big ( L \cap U' \big ) = L = \big ( L \cap W \big ) + \big ( L \cap W' \big )
\end{displaymath}
hold.
\end{lemma}
In the notations of the lemma above, we say that $W' \subseteq U'$ is a system of $L$-complements for $U \subseteq W$.
\par Now let 
\begin{displaymath}
0 = V_{0} < V_{1} < \cdots < V_{n} = V
\end{displaymath}
be a chain of subspaces of $V$ with $n \geq 1$. Then we consider the vector space
\begin{displaymath}
V^{\star} = V_{n-1} \oplus \frac{V}{V_{1}}
\end{displaymath}
We call it the derived vector space. Also, we have the full-dimensional lattice
\begin{displaymath}
L^{\star} = \big ( L \cap V_{n-1} \big ) + \frac{L + V_{1}}{V_{1}}
\end{displaymath}
of $V^{\star}$, which we call the derived lattice, and the chain of subspaces
\begin{displaymath}
0 = V^{\star}_{0} < V^{\star}_{1} < \cdots < V^{\star}_{n-1} = V^{\star}
\end{displaymath}
of $V^{\star}$ where
\begin{displaymath}
V^{\star}_{i} = V_{i-1} \oplus \frac{V_{i+1}}{V_{1}}
\end{displaymath}
which we call the derived chain. Note that its length is $n-1$, which is strictly less than the length of the chain of $V$.
\par Now let $W_{n-1} \subseteq W_{1}$ be a system of $L$-complements to $V_{1} \subseteq V_{n-1}$, and let us denote by $V \stackrel{\pi}{\longrightarrow} V_{1}$ the projection of $V$ onto $V_{1}$ along $W_{1}$. By $\mathrm{Lin}(W_{n-1},
V_1)$ we denote the space of all linear maps $W_{n-1}\to V_1$.
Then we consider the map
\begin{equation}\label{eq:epsilon}
\epsilon : \textrm{End}(V) \longrightarrow \textrm{Lin}(W_{n-1}, V_{1}) 
\qquad \varphi \mapsto \pi \circ \varphi_{| W_{n-1}}.
\end{equation}
We refer to it as the error map induced by the system $W_{n-1} \subseteq W_{1}$. Further we define
\begin{equation}\label{eq:Gamma}
\Gamma = \big \{ \gamma \in \textrm{Lin}(W_{n-1}, V_{1}) \ \vline \ \gamma(L \cap W_{n-1}) \subseteq L \cap V_{1} \big \}.
\end{equation}
Since $L \cap W_{n-1}$ is full-dimensional in $W_{n-1}$, $\Gamma$ is a lattice 
of $\textrm{Lin}(W_{n-1}, V_{1})$, and it is full-dimensional since 
$L \cap V_{1}$ is full-dimensional in $V_{1}$. We refer to it as the lattice 
induced by the system $W_{n-1} \subseteq W_{1}$.\par 
Now let $G$ be a unipotent algebraic group defined over $\mathbb{Q}$ acting 
faithfully on $V$, and suppose that
\begin{displaymath}
0 = V_{0} < V_{1} < \cdots < V_{n} = V
\end{displaymath}
is a flag of $V$ with respect to $G$, i.e., for all $v\in V_i$ we have 
$gv\equiv v \bmod V_{i-1}$ for all $g\in G$.
We consider the subgroup
\begin{displaymath}
G_{L} = \big \{ g \in G(\mathbb{Q}) \ \vline \ g L = L \big \}
\end{displaymath}
of $G(\mathbb{Q})$. We want to find a finite set of generators of this group.
\par
Since $V_{1}$ and $V_{n-1}$ are $G$-stable subspaces of $V$, $G$ acts on both 
$V_{n-1}$ and $\frac{V}{V_{1}}$, hence on their direct sum, that is to say, on 
the derived vector space. We refer to the action of $G$ on $V^{\star}$ as the 
derived action. Further, we denote by $N$ its kernel, which is of course a 
unipotent algebraic group over $\mathbb{Q}$ acting faithfully on $V$, by $Q$ 
its image, which is a unipotent algebraic group over $\mathbb{Q}$ acting 
faithfully on $V^{\star}$, and by $\pi$ the projection of $G$ onto $Q$. 
The following lemma is well-known; it follows directly from the commutativity
of diagram \eqref{eq:explog} in Section \ref{sect:basis}.

\begin{lemma}\label{lem:proj_surj}
The projection $\pi$ maps $G(\Q)$ surjectively onto $Q(\Q)$. 
\end{lemma}

Also we set
\begin{displaymath}
N_{L} = \big \{ g \in N(\mathbb{Q}) \ \vline \ g L = L \big \}
\end{displaymath}
and
\begin{displaymath}
Q_{L^{\star}} = \big \{ q \in Q(\mathbb{Q}) \ \vline \ q L^{\star} = L^{\star} \big \}.
\end{displaymath}
Of course the derived chain is a flag for $V^{\star}$ with respect to the action of $Q$. Since $G$ acts faithfully on $V$ we can regard elements in $G(\mathbb{Q})$ as automorphisms of $V$. The same consideration applies to elements of $N(\mathbb{Q})$. So we can apply the map $\epsilon$ to the elements of these groups.
\begin{proposition} \label{action}
For every $g \in G(\mathbb{Q})$ and $h \in N(\mathbb{Q})$ we have
\begin{displaymath}
\epsilon(g \cdot h) = \epsilon(g) + \epsilon(h),
\end{displaymath}
where $\epsilon$ is as in (\ref{eq:epsilon}).
\end{proposition}
\begin{proof}
Let $v \in W_{n-1}$. Since $h$ is an automorphism of $V$ acting as the identity on $\frac{V}{V_{1}}$,
\begin{displaymath}
h(v) - v \in V_{1}.
\end{displaymath}
Further, $g$ is an automorphism of $V$ acting as the identity on $V_{1}$, hence
\begin{displaymath}
g(h(v)- v) = h(v) - v.
\end{displaymath}
Since $W_{n-1} \subseteq W_{1}$, we have $\pi(v) = 0$. Hence applying $\pi$ to both sides of the previous identity and using linearity we obtain
\begin{displaymath}
\pi \circ g \circ h (v) = \pi \circ g (v) + \pi \circ h (v)
\end{displaymath}
hence the thesis.
\end{proof}
Of course, $N(\mathbb{Q})$ acts on $G(\mathbb{Q})$ by multiplication on the right. Once we endow $\textrm{Lin}(W_{n-1}, V_{1})$ with the obvious group structure given by addition, the previous proposition implies that the restriction of $\epsilon$ to $N(\mathbb{Q})$ is a group morphism. Hence $N(\mathbb{Q})$ acts on $\textrm{Lin}(W_{n-1}, V_{1})$ by
\begin{displaymath}
\textrm{Lin}(W_{n-1}, V_{1}) \times N(\mathbb{Q}) \rightarrow \textrm{Lin}(W_{n-1}, V_{1}) \qquad (x, h) \mapsto x + \epsilon(h).
\end{displaymath}
With these observations, we can restate the previous proposition saying that the restriction of $\epsilon$ to $G(\mathbb{Q})$ is a morphism of $N(\mathbb{Q})$-sets.
\par Now let us denote by $W$ the image of $N(\mathbb{Q})$ through $\epsilon$. Then Proposition \ref{action} gives us the
\begin{corollary} \label{commdiag}
There exists a unique map $\hat \epsilon : Q(\mathbb{Q}) \longrightarrow \frac{\mathrm{Lin}(W_{n-1}, V_{1})}{W}$ such that the diagram
\begin{displaymath}
\xymatrix{
G(\mathbb{Q}) \ar[r]^{\epsilon} \ar[d]^{\pi} & \mathrm{Lin}(W_{n-1}, V_{1}) \ar[d] \\
Q(\mathbb{Q}) \ar[r]^{\hat \epsilon} & \frac{\mathrm{Lin}(W_{n-1}, V_{1})}{W}
} 
\end{displaymath}
is commutative.
\end{corollary}
\begin{proof}
By Lemma \ref{lem:proj_surj}, the map $\pi : G(\Q) \rightarrow Q(\Q)$ 
is surjective. 
Also, if $g, g' \in G(\mathbb{Q})$ are such that $\pi(g) = \pi(g')$ then $g^{-1} g' \in N(\mathbb{Q})$, hence due to Proposition \ref{action} we obtain 
\begin{displaymath}
\epsilon(g') = \epsilon(g g^{-1} g') = \epsilon(g) + \epsilon(g^{-1} g')
\end{displaymath}
thus
\begin{displaymath}
\epsilon(g) + W = \epsilon(g') + W.
\end{displaymath}
These two facts show that the function
\begin{displaymath}
\hat \epsilon : Q(\mathbb{Q}) \longrightarrow \frac{\mathrm{Lin}(W_{n-1}, V_{1})}{W} \qquad q \mapsto \epsilon(g) + W
\end{displaymath}
where $g$ is any element of $G(\mathbb{Q})$ such that $\pi(g) = q$, is well defined and, of course, it makes the diagram above commutative. If $\hat \epsilon'$ is another such a function, then
\begin{displaymath}
\hat \epsilon \circ \pi = \hat \epsilon' \circ \pi
\end{displaymath}
hence, by surjectivity of $\pi$ it follows that 
$\hat \epsilon' = \hat \epsilon$.
\end{proof}
Also we set
\begin{displaymath}
G_{L^{\star}} = \big \{ g \in G(\mathbb{Q}) \ \vline \ g L^{\star} = L^{\star} \big \}.
\end{displaymath}
In other words, an element $g \in G(\mathbb{Q})$ lies in $G_{L^{\star}}$ if and only if $\pi(g)$ lies in $Q_{L^{\star}}$. Of course, $G_{L^{\star}}$ contains both $G_{L}$ and $N(\mathbb{Q})$.
\begin{lemma} \label{epsgamma}
Let $g \in G_{L^{\star}}$. Then  $g \in G_{L}$ if and only if $\epsilon(g) \in \Gamma$.
\end{lemma}
\begin{proof}
If $g \in G_{L}$ then $g(L \cap W_{n-1}) \subseteq L$, hence
\begin{displaymath}
\pi \circ g (L \cap W_{n-1}) \subseteq \pi (L) = \pi \big ( (L \cap V_{1}) \oplus (L \cap W_{1}) \big ) = L \cap V_{1}
\end{displaymath}
hence $\epsilon(g) \in \Gamma$.
\par Now let $g \in G_{L^{\star}}$ such that $\epsilon(g) \in \Gamma$. Then $g$ is an automorphism of $V$ fixing both $L + V_{1}$ and $L \cap V_{n-1}$, and such that $\pi \circ g$ sends $L \cap W_{n-1}$ in $L \cap V_{1}$. In particular, since the preimage of $L \cap V_{1}$ through $\pi$ is $W_{1} + (L \cap V_{1})$, we have that $g(L \cap W_{n-1}) \subseteq W_{1} + (L \cap V_{1})$. Further, since $L \cap W_{n-1} \subseteq L + V_{1}$ and $g$ fixes $L + V_{1}$, we have that $g(L \cap W_{n-1}) \subseteq L + V_{1}$. Hence
\begin{displaymath}
g(L \cap W_{n-1}) \subseteq \big ( L + V_{1} \big ) \cap \big ( W_{1} + (L \cap V_{1}) \big )
\end{displaymath}
Since $L = (L \cap V_{1}) + (L \cap W_{1})$, applying Dedekind's modular law we obtain equality
\begin{displaymath}
L = \big ( L + V_{1} \big ) \cap \big ( W_{1} + (L \cap V_{1}) \big )
\end{displaymath}
hence $g(L \cap W_{n-1}) \subseteq L$. Since $L = (L \cap V_{n-1}) + (L \cap W_{n-1})$ and $g$ fixes $L \cap V_{n-1}$, this shows that $g(L) \subseteq L$. Now let $l \in L$. Since $g$ fixes $V_{1} + L$, there exist $l_{1} \in V_{1}$ and $l_{2} \in L$ such that $g(l_{1}) + g(l_{2}) = l$. Since $g(L) \subseteq L$, we have that $g(l_{2}) \in L$, hence in particular that $g(l_{1}) = l - g(l_{2}) \in L$. Since $g$ fixes $V_{1}$, we also have that $g(l_{1}) \in V_{1}$, hence $g(l_{1}) \in L \cap V_{1}$. Since $V_{1} \subseteq V_{n-1}$, also $g(l_{1}) \in L \cap V_{n-1}$ holds. Since $g$ fixes $L \cap V_{n-1}$, we obtain that $l_{1} \in L \cap V_{n-1}$, hence $l_{1} + l_{2} \in L$, hence $L \subseteq g(L)$. So $g \in G_{L}$.
\end{proof}

\begin{proposition} \label{morphism}
The map given by the chain
\begin{displaymath}
G_{L^{\star}} \stackrel{\epsilon}{\longrightarrow} \mathrm{Lin}(W_{n-1}, V_{1}) \longrightarrow \frac{\mathrm{Lin}(W_{n-1}, V_{1})}{\Gamma}
\end{displaymath}
is a group morphism with kernel $G_{L}$.
\end{proposition}
\begin{proof}
Let $f, g \in G_{L^{\star}}$. Then they are automorphisms of $V$ acting as the identity on $V_{1}$, fixing $V_{n-1}$ and acting as the identity on $\frac{V}{V_{n-1}}$. Further, they fix $L \cap V_{n-1}$ and $L + V_{1}$. Now let $l \in L$. Of course, $g(l) - l \in V_{n-1} \cap (L + V_{1})$. Since $V_{1} \subseteq V_{n-1}$, applying Dedekind's modular law we obtain $V_{n-1} \cap (L + V_{1}) = V_{1} + (L \cap V_{n-1})$, which shows that
\begin{displaymath}
f \big ( g(l)- l \big ) - \big ( g(l) - l \big ) \in L \cap V_{n-1}
\end{displaymath}
hence
\begin{displaymath}
f \circ g (l) - f(l) - g(l) \in L.
\end{displaymath}
Since $\pi(L) = L \cap V_{1}$, we finally obtain
\begin{displaymath}
\pi \circ f \circ g (l) - \pi \circ f (l) - \pi \circ g (l) \in L \cap V_{1}
\end{displaymath}
which shows that the map is a group morphism. By Lemma \ref{epsgamma}, its kernel is $G_{L}$.
\end{proof}

\begin{proposition} \label{kerpsi}
Let $\hat\epsilon$ be as in Corollary \ref{commdiag}.
The map $\Psi$ given by the chain
\begin{displaymath}
Q_{L^{\star}} \stackrel{\hat \epsilon}{\longrightarrow} \frac{\mathrm{Lin}(W_{n-1}, V_{1})}{W} \longrightarrow \frac{\mathrm{Lin}(W_{n-1}, V_{1})}{W + \Gamma}
\end{displaymath}
is a group morphism. Its kernel is equal to the image of $G_{L}$ through $\pi$.
\end{proposition}
\begin{proof}
Since $\pi : G(\mathbb{Q}) \rightarrow Q(\mathbb{Q})$ is surjective and $g \in G(\mathbb{Q})$ is in $G_{L^{\star}}$ if and only if $\pi(g) \in Q_{L^{\star}}$, we obtain a surjective map $\pi : G_{L^{\star}} \rightarrow Q_{L^{\star}}$. By commutativity of the diagram in Corollary \ref{commdiag}, we also obtain the commutative diagram
\begin{displaymath}
\xymatrix{
G_{L^{\star}} \ar[r]^{\epsilon} \ar[d]^{\pi} & \textrm{Lin}(W_{n-1}, V_{1}) \ar[r] & \frac{\textrm{Lin}(W_{n-1}, V_{1})}{\Gamma} \ar[d] \\
Q_{L^{\star}} \ar[rr]^{\Psi} &  & \frac{\textrm{Lin}(W_{n-1}, V_{1})}{\Gamma + W}
}
\end{displaymath}
By Proposition \ref{morphism}, the top row is a group morphism, hence also the bottom row is. Again by Proposition \ref{morphism}, $G_{L}$ is the kernel of the top row. Hence, since the diagram above is commutative, the image of $G_{L}$ through $\pi$ lies in the kernel of $\Psi$. Now let $q$ be in the kernel of $\Psi$, and let $g \in G_{L^{\star}}$ be a preimage of $q$ through $\pi$. By commutativity of the diagram above, we have
\begin{displaymath}
\epsilon(g) \in \Gamma + W
\end{displaymath}
Now let $w \in W$ such that $\epsilon(g) + w \in \Gamma$. Since $W$ is the image of $N(\mathbb{Q})$ through $\epsilon$, there exists $h \in N(\mathbb{Q})$ such that $\epsilon(h) = w$. Hence by Proposition \ref{action} we have
\begin{displaymath}
\epsilon(g \cdot h) = \epsilon(g) + \epsilon(h) \in \Gamma
\end{displaymath}
thus by Lemma \ref{epsgamma} we obtain that $g \cdot h \in G_{L}$. Of course, $\pi(g \cdot h) = q$.
\end{proof}

\section{An example}\label{sect:exa}

Let us consider
\begin{displaymath}
G = \left \{ \left ( \begin{array}{cccc}
1 & 0 & a & b \\
0 & 1 & c & \frac{1}{2}c^{2} \\
0 & 0 & 1 & c \\
0 & 0 & 0 & 1
\end{array} \right ) \in \textrm{GL}_{4}(\mathbb{C}) \ \textrm{such that $a, b, c \in \mathbb{C}$} \right \}.
\end{displaymath}
It is easy to check that $G$ is an algebraic subgroup of $\textrm{GL}_{4}(\mathbb{C})$ defined over $\mathbb{Q}$. 
%For example, it is the Zariski-closed of $\textrm{GL}_{n}(\mathbb{C})$ 
%associated to the ideal
%\begin{displaymath}
%\left ( \begin{array}{c} 
%\{ x_{i,j} \ \vline \ 1 \leq j < i \leq 4 \} \ \cup \ \{ x_{i,i} - 1 \ \vline \ i = 1, \dots, 4 \} \ \cup \\
%\cup \ \{ x_{2,1}, \ x_{2,3} - x_{3,4}, \ x_{2,3} x_{3,4} - 2 x_{2,4} \}
%\end{array} \right )
%\end{displaymath}
Since it is contained in the set of upper-unitriangular matrices of 
$\textrm{GL}_{4}(\mathbb{C})$, it is unipotent. In this example it is rather
straightforward to find a set of generators of $G(\Z)$ directly. However, in 
this section
we illustrate the results of the previous section by showing how they 
help us finding a finite set of generators for $G(\Z)$.
\par 
$G$ acts faithfully on $\mathbb{Q}^{4}$ by matrix-vector multiplication; 
also, $\mathbb{Z}^{4}$ is a full-dimensional lattice of $\mathbb{Q}^{4}$. 
Thus we can consider the subgroup $G_{\mathbb{Z}^{4}}$ of $G$, and it is easily 
seen that
\begin{displaymath}
G(\Z) = G_{\mathbb{Z}^{4}}.
\end{displaymath}
The chain of subspaces
\begin{displaymath}
0 = V_{0} < \langle e_{1}, e_{2} \rangle = V_{1} < \langle e_{1}, e_{2}, e_{3} \rangle = V_{2} < V_{3} = \mathbb{Q}^{4}
\end{displaymath}
where $e_{1}$, $e_{2}$, $e_{3}$ and $e_{4}$ are the standard basis of 
$\mathbb{Q}^{4}$, is a flag of $V$ with respect to the action of $G$. The derived vector space has basis given by $(e_{1}, 0)$, $(e_{2}, 0)$, $(e_{3}, 0)$, $(0, e_{3} + V_{1})$ and $(0, e_{4} + V_{1})$; through this basis we can identify it with $\mathbb{Q}^{5}$, and under this identification the derived lattice corresponds to $\mathbb{Z}^{5}$. The kernel of the derived action is
\begin{displaymath}
N = \left \{ \left ( \begin{array}{cccc}
1 & 0 & 0 & b \\
0 & 1 & 0 & 0 \\
0 & 0 & 1 & 0 \\
0 & 0 & 0 & 1
\end{array} \right ) \in \textrm{GL}_{4}(\mathbb{C}) \ \textrm{such that $b \in \mathbb{C}$} \right \}.
\end{displaymath}
In particular,
\begin{displaymath}
N_{\mathbb{Z}^{4}} = \left \{ \left ( \begin{array}{cccc}
1 & 0 & 0 & b \\
0 & 1 & 0 & 0 \\
0 & 0 & 1 & 0 \\
0 & 0 & 0 & 1
\end{array} \right ) \in \textrm{GL}_{4}(\mathbb{C}) \ \textrm{such that $b \in \mathbb{Z}$} \right \}
\end{displaymath}
and it is straighforward to check that it is an infinite cyclic group with generator
\begin{displaymath}
n = \left ( \begin{array}{cccc}
1 & 0 & 0 & 1 \\
0 & 1 & 0 & 0 \\
0 & 0 & 1 & 0 \\
0 & 0 & 0 & 1
\end{array} \right ).
\end{displaymath}
The subgroup
\begin{displaymath}
Q = \left \{ \left ( \begin{array}{ccc|cc}
1 & 0 & a & 0 & 0 \\
0 & 1 & c & 0 & 0 \\
0 & 0 & 1 & 0 & 0 \\
\hline
0 & 0 & 0 & 1 & c \\
0 & 0 & 0 & 0 & 1
\end{array} \right ) \in \textrm{GL}_{5}(\mathbb{C}) \ \textrm{such that $a, c \in \mathbb{C}$} \right \}
\end{displaymath}
of $\mathrm{GL}_{5}(\mathbb{C})$ is the image of the derived action, the projection of $G$ onto $Q$ being
\begin{displaymath}
\pi : G \rightarrow Q \quad \left ( \begin{array}{cccc} 
1 & 0 & a & b \\
0 & 1 & c & \frac{1}{2}c^{2} \\
0 & 0 & 1 & c \\
0 & 0 & 0 & 1
\end{array} \right ) \mapsto \left ( \begin{array}{ccc|cc}
1 & 0 & a & 0 & 0 \\
0 & 1 & c & 0 & 0 \\
0 & 0 & 1 & 0 & 0 \\
\hline
0 & 0 & 0 & 1 & c \\
0 & 0 & 0 & 0 & 1
\end{array} \right ). 
\end{displaymath}
%Under the previous correspondence between the derived vector space and 
%$\mathbb{Q}^{5}$ the action of $Q$ on the derived vector space corresponds to 
%the action of $Q$ on $\mathbb{Q}^{5}$ given by multiplication matrix-vector. 
Now $Q_{\mathbb{Z}^{5}}$ is a torsion free abelian group of rank 2 with basis 
given by
\begin{displaymath}
q_{1} = \left ( \begin{array}{ccc|cc}
1 & 0 & 1 & 0 & 0 \\
0 & 1 & 0 & 0 & 0 \\
0 & 0 & 1 & 0 & 0 \\
\hline
0 & 0 & 0 & 1 & 0 \\
0 & 0 & 0 & 0 & 1
\end{array} \right ) \qquad q_{2} = \left ( \begin{array}{ccc|cc}
1 & 0 & 0 & 0 & 0 \\
0 & 1 & 1 & 0 & 0 \\
0 & 0 & 1 & 0 & 0 \\
\hline
0 & 0 & 0 & 1 & 1 \\
0 & 0 & 0 & 0 & 1
\end{array} \right ).
\end{displaymath}
A system of $\mathbb{Z}^{4}$-complements for $V_{1} \subseteq V_{2}$ is given by $W_{2} \subseteq W_{1}$ where
\begin{displaymath}
W_{2} = \langle e_{4} \rangle \qquad W_{1} = \langle e_{3}, e_{4} \rangle
\end{displaymath}
Using the basis $e_{4}$ for $W_{2}$ and the basis $e_{1}$, $e_{2}$ for $V_{1}$ we can identify $\textrm{Lin}(W_{2}, V_{1})$ with $\textrm{M}_{2 \times 1}(\mathbb{Q})$; under this identification, the induced lattice $\Gamma$ corresponds to $\textrm{M}_{2 \times 1}(\mathbb{Z})$. Further, using the standard basis of $\mathbb{Q}^{4}$ we can identify $\textrm{End}(\mathbb{Q}^{4})$ with $\textrm{M}_{4 \times 4}(\mathbb{Q})$. In this way, the error map is
\begin{displaymath}
\epsilon : \textrm{M}_{4 \times 4}(\mathbb{Q}) \rightarrow  \textrm{M}_{2 \times 1}(\mathbb{Q}) \qquad \left ( \begin{array}{cccc}
a_{1,1} & a_{1,2} & a_{1,3} & a_{1,4} \\
a_{2,1} & a_{2,2} & a_{2,3} & a_{2,4} \\
a_{3,1} & a_{3,2} & a_{3,3} & a_{3,4} \\
a_{4,1} & a_{4,2} & a_{4,3} & a_{4,4}
\end{array} \right ) \mapsto \left ( \begin{array}{c}
a_{1,4} \\
a_{2,4}
\end{array} \right ).
\end{displaymath}
The image of the rational points of $N$ through $\epsilon$ is the subspace of 
$\mathrm{M}_{2 \times 1}(\mathbb{Q})$ generated by $\binom{1}{0}$. So in the 
notation of Section \ref{thresults} we have $W=\langle \binom{1}{0} \rangle$
and $\Gamma = \mathrm{M}_{2\times 1}(\Z)$. Furthermore, the map $\Psi$ from
Proposition \ref{kerpsi} goes from the rational points of $Q$ to 
$\frac{\mathrm{M}_{2 \times 1}(\mathbb{Q})}{W+\Gamma}$.\par 
Now we need two matrices $g_{1}$ and $g_{2}$ in $G(\mathbb{Q})$ 
whose images through $\pi$ are $q_{1}$ and $q_{2}$, respectively. Their 
existence is guaranteed by the surjectivity of $\pi$. For example, we can take
\begin{displaymath}
g_{1} = \left ( \begin{array}{cccc}
1 & 0 & 1 & 0 \\
0 & 1 & 0 & 0 \\
0 & 0 & 1 & 0 \\
0 & 0 & 0 & 1
\end{array} \right ) \qquad g_{2} = \left ( \begin{array}{cccc}
1 & 0 & 0 & 0 \\
0 & 1 & 1 & \frac{1}{2} \\
0 & 0 & 1 & 1 \\
0 & 0 & 0 & 1
\end{array} \right ).
\end{displaymath}
This shows in particular that $g_{1}$ and $g_{2}$ are in $G_{\mathbb{Z}^{5}}$. 
Also, using commutativity of the diagram in Corollary \ref{commdiag}, we have 
that
\begin{displaymath}
\Psi(q_{1}) = 0 +W+\Gamma \qquad \Psi(q_{2}) = \left ( \begin{array}{c}
0 \\
\frac{1}{2}
\end{array} \right ) + W+\Gamma.
\end{displaymath}
Now the kernel of $\Psi$ is generated by $q_{1}$ and 
$q_{2}^{2}$. Therefore, by Proposition \ref{kerpsi}, $\pi(G_{\Z^4})$ is 
generated by $q_1$ and $q_2^2$. Their preimages, $g_1$ and $g_2^2$ are
only guaranteed to lie in $G_{\Z^5}$; however, here we see that they are
alreay in $G_{\Z^4}$. Now the kernel of $\pi$ restricted to $G_{\Z^4}$ is
$N_{\Z^4}$. So $g_{1} N_{\mathbb{Z}^{4}}$ and $g_{2}^{2} N_{\mathbb{Z}^{4}}$ generate
$\frac{G_{\mathbb{Z}^{4}}}{N_{\mathbb{Z}^{4}}}$. We conclude that $n,g_1,g_2^2$ 
generate $G_{\Z^4}$.\par 
Of course we could have made a different choice for the preimages of 
$q_1$ and $q_2$. For example, we could have taken
\begin{displaymath}
g'_{1} = \left ( \begin{array}{cccc}
1 & 0 & 1 & \frac{1}{2} \\
0 & 1 & 0 & 0 \\
0 & 0 & 1 & 0 \\
0 & 0 & 0 & 1
\end{array} \right )
\end{displaymath}
instead of $g_1$. 
Then $g'_{1}$ is in $G_{\mathbb{Z}^{5}}$ but it is not in $G_{\mathbb{Z}^{4}}$. 
However, since $q_{1}$ is in the kernel of $\Psi$, we get $\hat\epsilon(q_1)
\in W+\Gamma$ and by the commutativity of the diagram in Corollary 
\ref{commdiag}, $\epsilon(g_{1}') \in W+\Gamma$. This can of course also
be checked directly as 
$$\epsilon(g_1') = \begin{pmatrix} \tfrac{1}{2} \\ 0 \end{pmatrix}.$$
Now we note that
\begin{displaymath}
n' = \left ( \begin{array}{cccc}
1 & 0 & 0 & \frac{1}{2} \\
0 & 1 & 0 & 0 \\
0 & 0 & 1 & 0 \\
0 & 0 & 0 & 1
\end{array} \right ) \in N(\mathbb{Q}) \text{ with }
\epsilon(n') = \left ( \begin{array}{c}
 \frac{1}{2} \\
0
\end{array} \right ).
\end{displaymath}
The existence of such an $n'$ is guaranteed by the fact that
$W$ is the image of the rational points of $N$ through $\epsilon$. 
Since $N$ is in the kernel of $\pi$ we have 
$\pi(g'_{1} \cdot (n')^{-1}) = \pi(g_1')=q_{1}$. So we can work with 
$g_1'(n')^{-1}$ as preimage of $q_1$. This choice of preimage works for us
as, due to Proposition \ref{action} we have that
\begin{displaymath}
\epsilon(g'_{1} \cdot (n')^{-1}) = \epsilon(g'_{1}) - \epsilon(n') = 
\left ( \begin{array}{c}
\frac{1}{2} \\
0
\end{array} \right ) - \left ( \begin{array}{c}
\frac{1}{2} \\
0
\end{array} \right ) = \left ( \begin{array}{c}
0 \\
0
\end{array} \right ) \in\Gamma.
\end{displaymath}
Hence $g'_{1} (n')^{-1} \in G_{\mathbb{Z}^{4}}$ and we can apply the same 
considerations as above to prove that $n$, $g'_{1} (n')^{-1}$ and $g_{2}^{2}$ 
are a generating set for $G_{\mathbb{Z}^{4}}$. As a matter of coincidence, we 
note that
\begin{displaymath}
g'_{1} \cdot (n')^{-1} = g_{1}.
\end{displaymath}
Finally we note that $n$ commutes with both $g_{1}$ and $g_{2}^{2}$, hence the chain of subgroups
\begin{displaymath}
1 < \langle n \rangle < \langle n, g_{1} \rangle < \langle n, g_{1}, g_{2}^{2} \rangle = G_{\mathbb{Z}^{4}}
\end{displaymath}
is a central series for $G_{\mathbb{Z}^{4}}$ with infinite cyclic factors.

\section{The Lie algebra connection}\label{sect:basis}

Let $G\subset \GL_m(\C)$ be a unipotent algebraic group. As illustrated in
Section \ref{sect:exa}, the results in Section \ref{thresults} in principle 
yield an algorithm for finding a finite set of generators of $G(\Z)$.
However, to make it work efficiently in practice we rather work with
the Lie algebra of $G$ than with $G$ itself. In this section we describe
the main reults that we need for that. \par  
First we review some standard facts on the Lie algebra of an 
algebraic group; for more details we refer to \cite{borel}, 
\cite{chevalley}, \cite{tauvelyu},\cite{waterhouse}.\par
As customary we denote the Lie algebras of the algebraic groups $G,H,\ldots$
by $\mathfrak{g}, \mathfrak{h},\ldots$. Let $G$ be a unipotent algebraic
group defined over $\Q$, acting on a vector space $V$. Then $\g$ also
acts on $V$. Furthermore, $G$ is connected, and hence a subspace $U\subset V$
is $G$-stable if and only if it is $\g$-stable. If this is the case then
we get a $G$-action and a $\g$-action on $U$, and those are compatible, in
the sense that the corresponding $\g$-representation is the differential
of the $G$-representation. Similarly we get compatible $G$- and $\g$-actions
on quotients and direct sums of modules. \par
An important role in our algorithm is played by the exponential mapping.
For a nilpotent $x\in \gl_m(\C)$ we set 
$$ \exp(x) = \sum_{i=0}^{n-1} \frac{x^i}{i!},$$
and for a unipotent $u\in \GL_m(\C)$ 
$$\log(u) = \sum_{i=1}^{n-1} (-1)^{i-1}\frac{(u-1)^i}{i}.$$
Since $G$ is defined over $\Q$ we have that $\g\subset \gl_m(\C)$ has a basis 
such that all elements have coefficients in $\Q$. The $\Q$-span of such a
basis is denoted $\g_\Q$. Then it is well-known that the maps
$\exp : \g_\Q\to G(\Q)$ and $\log : G(\Q)\to \g_\Q$ are mutually inverse. 
In particular, when we work with the Lie algebra $\g_\Q$ we keep control 
over the elements of $G(\Q)$ by these mappings. \par
Now let $W$ be another finite dimensional vector space over $\mathbb{Q}$, 
$H$ a unipotent algebraic subgroup of $\GL(W)$, and
$\varphi : G \rightarrow H$ a morphism of algebraic groups. 
Then we have the diagram
\begin{equation}\label{eq:explog}
\xymatrix{
G(\mathbb{Q}) \ar@/^/[rr]^{\textrm{log}} \ar[d]^{\varphi} && \mathfrak{g} \ar@/^/[ll]^{\textrm{exp}} \ar[d]^{\textrm{d} \varphi} \\
H(\mathbb{Q}) \ar@/^/[rr]^{\textrm{log}} && \mathfrak{h} \ar@/^/[ll]^{\textrm{exp}}
}
\end{equation}
and it turns out that it is commutative, i.e., $\exp(\mathrm{d}\varphi(x))
=\varphi (\exp(x))$ for all $x\in \mathfrak{g}$ (cf. \cite{chevalley},
Chapter V, \S 4, Proposition 15). 

Now we return to the setting of Section \ref{thresults}.
A sequence of subspaces 
$$0=V_0 < V_1 <\cdots <V_n=V$$
is a flag for the action of $G$ if and only if it is a flag for the
action of $\g$. (The latter means that $\g\cdot V_i \subset V_{i-1}$ for
$i>0$.) In particular, $\g$ acts on the derived vector space $V^\star$,
and the corresponding representation of $\g$ is the differential of the
representation of $G$ on $V^\star$. In particular this means that 
$\mathfrak{n}$, which is the Lie algebra of $N$, is the kernel of 
$\mathrm{d}\pi$ and $\mathfrak{q} = \mathrm{d}\pi (\g)$, where $\mathfrak{q}$
is the Lie algebra of $Q$. 

\begin{proposition} \label{central}
$\mathfrak{n}$ is central in $\mathfrak{g}$.
\end{proposition}
\begin{proof}
Let $x \in \mathfrak{g}$ and $y \in \mathfrak{n}$. Then $x$ is an endomorphism 
of $V$ such that $x.V \subseteq V_{n-1}$ and $x.V_{1} = 0$, and $y$ is an 
endomorphism of $V$ such that $y.V \subseteq V_{1}$ and $y.V_{n-1} = 0$. Now 
let $v \in V$. Then $y.v \in V_{1}$, hence $x.y.v = 0$. Also, $y.v \in V_{1}$, 
hence $x.y.v = 0$. Thus $[x, y].v = 0$. Since $\mathfrak{g}$ acts faithfully 
on $V$, $[x, y] = 0$.
\end{proof}

Now let $L$ be a full-dimensional lattice of $V$, $W_{n-1} \subseteq W_{1}$ a 
system of $L$-complements to $V_{1} \subseteq V_{n-1}$, and let us denote by 
$\epsilon$ and by $\Gamma$ the induced error map and the induced lattice, 
respectively. Since $\mathfrak{g}$ acts faithfully on $V$, we can regard its 
elements as endomorphisms of $V$. The same consideration applies to 
$\mathfrak{n}$. So we can consider the restriction of the map $\epsilon$
to $\g$ and $\mathfrak{n}$.

\begin{proposition}\label{prop:ninj}
The restriction of the induced map $\epsilon$ to $\mathfrak{n}$ is injective.
\end{proposition}
\begin{proof}
Let $x \in \mathfrak{n}$ such that $\epsilon(x) = 0$. Since $x(V) \subseteq V_{1}$, for every $v \in V$ we have
\begin{displaymath}
x(v) = \pi \circ x (v)
\end{displaymath}
Since $\epsilon(x) = 0$, $\pi \circ x (W_{n-1}) = 0$. Hence $x(W_{n-1}) = 0$. Since $x(V_{n-1}) = 0$, we finally have $x(V) = 0$.
\end{proof}
Further, we define
\begin{displaymath}
\mathfrak{n}_{L}  = \big \{ x \in \mathfrak{n} \ \vline \ \epsilon(x) \in \Gamma \big \}.
\end{displaymath}
By the previous proposition, it is a full-dimensional lattice of $\mathfrak{n}$. Since $\mathfrak{n}$ acts faithfully on $V$ and $V$ admits a flag with respect to this action, then $\mathfrak{n}$, regarded as a Lie subalgebra of $\mathfrak{gl}(V)$, consists of nilpotent endomorphisms. Hence we can consider the diagram
\begin{displaymath}
\xymatrix{
\mathfrak{n} \ar[r]^{\textrm{exp}} \ar[rd]_{\epsilon} & \GL(V) \ar[d]^{\epsilon} \\
& \textrm{Lin}(W_{n-1}, V_{1})
}
\end{displaymath}
\begin{proposition} \label{commtriang}
The diagram above is commutative.
\end{proposition}
\begin{proof}
Let us denote by $\textrm{id}_{V}$ the identity endomorphism of $V$. Since $W_{n-1} \subseteq W_{1}$, $\epsilon(\textrm{id}_{V}) = 0$. Now let $x \in \mathfrak{n}$. Then $x(V) \subseteq V_{1}$ and $x(V_{n-1}) = 0$, hence $x^{2}(V) = 0$, and
\begin{displaymath}
\textrm{exp} \ x = \textrm{id}_{V} + x
\end{displaymath}
Thus
\begin{displaymath}
\epsilon(\textrm{exp} \ x) = \epsilon(\textrm{id}_{V} + x) = \epsilon(\textrm{id}_{V}) + \epsilon(x) = \epsilon(x)
\end{displaymath}
\end{proof}

\section{$T$-groups}\label{sect:tgrps}

The groups $G(\Z)$ that we are after are finitely-generated nilpotent
and torsion free. Such groups are called $T$-groups in the literature
(cf. \cite{segal}). In this section we review some facts that we need on
$T$-groups. \par
It is known that any $T$-group admits a (proper normal) central 
series with infinite cyclic factors; the other way round, every group 
admitting such a series is clearly a $T$-group. Now let $G$ be a group, and 
let $g_{1}, \dots, g_{n}$ be a (ordered) set of elements of $G$. Then we can 
consider the chain of subgroups
\begin{displaymath}
G_{1} \geq G_{2} \geq \cdots \geq G_{n} \geq G_{n+1} = 1
\end{displaymath}
of $G$ where for every $i = 1, \dots, n$,
\begin{displaymath}
G_{i} = \langle g_{i}, \dots, g_{n} \rangle
\end{displaymath}
We call it the chain associated to $g_{1}, \dots, g_{n}$. Further, we say 
that $g_{1}, \dots, g_{n}$ is a $T$-sequence for $G$ if the associated chain 
is a proper central series for $G$ with infinite cyclic factors. It will be 
convenient to extend this terminology saying that the empty set is a 
$T$-sequence for the trivial group. Every $T$-group $G$ has a $T$-sequence, and
the length of a $T$-sequence is an invariant of the group, called the 
Hirsch-length of $G$. Furthermore, if $g_1,\ldots,g_n$ is a $T$-sequence 
for $G$ then every element $g\in G$ can be written $g = g_1^{e_1}\cdots
g_n^{e_n}$, where $e_i\in \Z$.

Now let $A$ be an abelian group, and let $a_{1}, \dots, a_{n} \in A$. Then we 
consider
\begin{displaymath}
L = \big \{ (e_{1}, \dots, e_{n}) \in \mathbb{Z}^{n} \ \vline \ a_{1}^{e_{1}} \cdots a_{n}^{e_{n}} = 1 \big \}
\end{displaymath}
which is of course a subgroup of $\mathbb{Z}^{n}$. We call it the relation 
lattice of $a_{1}, \dots, a_{n}$ in $A$. For every $e = (e_{1}, \dots, e_{n}) 
\in L$, $e \neq 0$, we can define its height as the minimum $j = 1, \dots, n$ 
such that $e_{j} \neq 0$, and its leading coefficient as the integer $e_{j}$. 
Now let $e^{(1)}, \dots, e^{(m)}$ be a basis of $L$, where for 
$j = 1, \dots, m$
\begin{displaymath}
e^{(j)} = ( e^{(j)}_{1}, \dots, e^{(j)}_{n}).
\end{displaymath}
We say that the basis is in Hermite normal form if the matrix
\begin{displaymath}
\left ( \begin{array}{ccc}
e^{(1)}_{1} & \cdots & e^{(1)}_{n} \\
\vdots & & \vdots \\
e^{(m)}_{1} & \cdots & e^{(m)}_{n}
\end{array} \right ) \in M_{m \times n}(\mathbb{Z})
\end{displaymath}
is. This means that there exists a sequence of integers
\begin{displaymath}
1 \leq i_{1} < \cdots < i_{m} \leq n
\end{displaymath}
such that
\begin{displaymath}
e^{(j)}_{i} = 0
\end{displaymath}
for all $j = 1, \dots, m$ and all $1 \leq i < i_{j}$, and that
\begin{displaymath}
0 \leq e^{(k)}_{i_{j}} < e^{(j)}_{i_{j}}
\end{displaymath}
for every $1 \leq k < j \leq m$. It is known that $L$ admits a unique 
basis in Hermite normal form.
We note that the height of any non-zero element of 
$L$ is any of the integers $i_{1}, \dots, i_{m}$. Further, if its height is 
$i_{j}$ for some $j = 1, \dots, m$, then its leading coefficient is a 
(non-zero) multiple of $e^{(j)}_{i_{j}}$. \par 
The next lemma is a somewhat stronger version of a result of Eick 
(cf. \cite{eick}, Lemma 3.19).

\begin{lemma} \label{lemma:Eick}
Let $G$ be a $T$-group, $A$ an abelian group, and let $\varphi : G \rightarrow 
A$ be a morphim of groups. Further, let $g_{1}, \dots, g_{n}$ be a $T$-sequence
for $G$, and let $e^{(1)}, \dots, e^{(m)}$ be the basis in Hermite normal form 
of the relation lattice of $\varphi(g_{1}), \dots, \varphi(g_{n})$ in $A$. 
For $j = 1, \dots, m$ set
\begin{displaymath}
k_{j} = g_{1}^{e^{(j)}_{1}} \cdots g_{n}^{e^{(j)}_{n}}.
\end{displaymath}
Then $k_1,\ldots,k_m$ 
is a $T$-sequence for the kernel of $\varphi$.
\end{lemma}
\begin{proof}
Let us denote by
\begin{displaymath}
K_{1} \geq K_{2} \geq \cdots \geq K_{m} \geq K_{m+1} = 1
\end{displaymath}
the chain of subgroups of $G$ associated to $k_{1}, \dots, k_{m}$. We want to 
show that it is a (proper normal) central series for $\mathrm{ker} \varphi$ 
with infinite cyclic factors. Since the basis $e^{(1)}, \dots, e^{(m)}$ is in 
Hermite normal form, we have the sequence
\begin{displaymath}
1 \leq i_{1} < i_{2} < \cdots < i_{m} \leq n
\end{displaymath}
defined as above. It is convenient to set $i_{0} = 0$ and $i_{m+1} = n + 1$. Also, let us denote by
\begin{displaymath}
G = G_{1} > G_{2} > \cdots > G_{n} > G_{n+1} = 1
\end{displaymath}
the proper central sequence with infinite cyclic factors for $G$ associated to $g_{1}, \dots, g_{n}$. Then it is enough to prove that for every $j = 1, \dots, m+1$ we have
\begin{displaymath}
\mathrm{ker} \varphi \cap G_{i_{j-1} + 1} = \cdots = \mathrm{ker} \varphi \cap G_{i_{j}-1} = \mathrm{ker} \varphi \cap G_{i_{j}} = K_{j}
\end{displaymath}
Indeed, suppose that the previous equalities hold. Then $K_{1} = \mathrm{ker} \varphi \cap G_{i_{0} + 1} = \mathrm{ker} \varphi \cap G_{1} = \mathrm{ker} \varphi \cap G = \mathrm{ker} \varphi$. Also, for every $j = 1, \dots, m+1$, $G_{i_{j}} \unlhd G$, hence $K_{j} = \mathrm{ker} \varphi \cap G_{i_{j}} \unlhd \mathrm{ker} \varphi \cap G = \mathrm{ker} \varphi$. This shows that the chain is a normal series for $\mathrm{ker} \varphi$. Further, for every $j = 1, \dots, m$, the map $\mathrm{ker} \varphi \rightarrow \frac{G}{G_{i_{j} + 1}}$ has kernel $\mathrm{ker} \varphi \cap G_{i_{j}+1} = K_{j+1}$. Hence it factors through a group monomorphism $\frac{\mathrm{ker} \varphi}{K_{j+1}} \rightarrow \frac{G}{G_{i_{j} + 1}}$. The image of $\frac{K_{j}}{K_{j+1}}$ through it is $\frac{K_{j} G_{i_{j}+1}}{G_{i_{j}+1}} = \frac{(\mathrm{ker} \varphi \cap G_{i_{j}})G_{i_{j}+1}}{G_{i_{j}+1}} = \frac{(\mathrm{ker} \varphi G_{i_{j}+1}) \cap G_{i_{j}}}{G_{i_{j}+1}} \leq \frac{G_{i_{j}}}{G_{i_{j}+1}}$, and the image of $k_{j} K_{j+1}$ is $g^{e_{i_{j}}^{(j)}}_{i_{j}} G_{i_{j}+1}$. This shows that the series is central and with infinite cyclic factors.
\par So we have to prove the previous equalities. It is clear that for every $j = 1, \dots, m+1$,
\begin{displaymath}
\mathrm{ker} \varphi \cap G_{i_{j-1}+1} \supseteq \cdots \supseteq \mathrm{ker} \varphi \cap G_{i_{j}-1} \supseteq \mathrm{ker} \varphi \cap G_{i_{j}} \supseteq K_{j}
\end{displaymath}
and it remains to prove the reverse inclusions. We proceed by induction on $j$. Let us consider the base case $j = m+1$. Then $K_{j}$ is trivial, and all we have to show is that for every $l = i_{m} + 1, \dots, n+1$, $\mathrm{ker} \varphi \cap G_{l}$ is trivial, too. Again, we proceed by induction on $l$. In the base case $l = n+1$, it is obviously true. Now let $l = i_{m} + 1, \dots, n$ and suppose that $\mathrm{ker} \varphi \cap G_{l+1}$ is trivial.  Let $g \in \mathrm{ker} \varphi \cap G_{l}$. Then $g = g_{l}^{e} h$ for some $e \in \mathbb{Z}$ and $h \in G_{l+1}$, and 
\begin{displaymath}
\varphi(g_{l})^{e} + \varphi(h) = 0
\end{displaymath}
Since $h \in \langle g_{l+1}, \dots, g_{n}\rangle$, then $\varphi(h) \in \langle \varphi(g_{l+1}), \dots, \varphi(g_{n}) \rangle$. Thus if $e \neq 0$, then there would exist an element in the relation lattice with height $l$, which is impossible. Hence $e = 0$, hence $g = h \in \mathrm{ker} \varphi \cap G_{l+1}$, hence $g = 1$ by the inductive hypothesis. This concludes the case $j = m+1$. Now let $j = 1, \dots, m$, and suppose that
\begin{displaymath}
\mathrm{ker} \varphi \cap G_{i_{j} + 1} = \cdots = \mathrm{ker} \varphi \cap G_{i_{j+1}-1} = \mathrm{ker} \varphi \cap G_{i_{j+1}} = K_{j+1}
\end{displaymath}
In this case we have to show that for every $l = i_{j-1}+1, \dots, i_{j}$,
\begin{displaymath}
\mathrm{ker} \varphi \cap G_{l} = K_{j}
\end{displaymath}
and again we proceed by induction on $l$. Let us just consider the base case $l = i_{j}$, the inductive step being similar to the one in the case $j = m+1$. Let $g \in \mathrm{ker} \varphi \cap G_{i_{j}}$. Then $g = g_{i_{j}}^{e} h$ for some $e \in \mathbb{Z}$ and some $h \in G_{i_{j}+1}$. If $e = 0$ then $g \in G_{i_{j}+1}$ and we conclude by inductive hypothesis that $g \in K_{j+1}$. Now let us suppose $e \neq 0$. Arguing as before, the relation lattice contains an element of height $i_{j}$ and leading coefficient $e$. Thus $e_{i_{j}}^{(j)}$ divides $e$. Let us denote by $f$ the quotient. Then $g G_{i_{j}+1} = k_{i}^{f} G_{i_{j}+1}$, hence by inductive hypothesis $g k_{j}^{-f} \in \mathrm{ker} \varphi \cap G_{i_{j}+1} = K_{j+1}$, hence finally $g \in K_{j}$.
\end{proof}

\section{Some algorithms for lattices}\label{sect:smith}

In this section we describe some algorithms that solve several problems 
related to lattices. We mainly work with matrices whose rows span lattices
or subspaces in $\Q^n$. We say that a matrix is integral if it has integer
entries.\par
A basic algorithm that we use is the Smith normal form:
given an $m\times n$ integral matrix $A$ this algorithm finds an $m\times
n$ integral matrix $S$, and integral unimodular square matrices $P$ and $Q$ 
with
\begin{enumerate}
\item $S$ is in Smith normal form (this means that there is an $r$ such
that $d_i = S(i,i)$ is positive for $1\leq i\leq r$, $S$ has no other 
nonzero entries, and $d_i$ divides $d_{i+1}$ for $1\leq i <r$),
\item $S=PAQ$.
\end{enumerate}
For details on this algorithm we refer to \cite{sims}, \S 8.3. One property
that we note is the following.

\begin{lemma}\label{lem:snf}
Let $q_1,\ldots,q_n$ denote the rows of $Q^{-1}$. They form a basis of $\Z^n$,
and $d_i$ is the smallest non-negative integer such that 
$d_iq_i$ lies in the span of the rows of $A$ for $1\leq i\leq r$. 
\end{lemma}

We also need an algorithm that appears to be well known: in 
the computer algebra system {\sc Magma} (\cite{magma}) it
is implemented under the name of ``saturation''. However, we have not been
able to find a reference for it in the literature. For this reason we sketch 
a solution here. Let $A$ be an $m\times n$-matrix with integer entries. Let
$V\subset \Q^n$ be the $\Q$-space spanned by the rows of $A$. The problem
is to find a $\Z$-basis for the lattice $\Z^n\cap V$. Without loss of
generality we assume that the rows of $A$ are linearly independent.
The key observation is the following. Let $B$ be an $m\times n$ integral
matrix whose rows span $V$. Then its rows span $\Z^n \cap V$
if and only if the Smith normal form of $B$ has diagonal entries that are
all equal to $1$. This follows from Lemma \ref{lem:snf}. This yields the
following algorithm.

\begin{algorithm}[Saturation]\label{alg:saturation}
\noindent{\bf Input:} an $m\times n$ integral matrix $A$ with linearly
independent rows.\\
{\bf Output:} an $m\times n$ integral matrix $B$ whose rows span $\Z^n\cap V$,
where $V\subset \Q^n$ is the $\Q$-space spanned by the rows of $A$.
\begin{enumerate}
\item Let $S$, $P$, $Q$ be the output of the Smith normal form algorithm with
input $A$.
\item Let $S'$ be the matrix obtained from $S$ by setting the diagonal entries
equal to $1$.
\item Return $B= P^{-1}S'Q^{-1}$.
\end{enumerate}
\end{algorithm}

\begin{algorithm}[Intersection of lattice and subspace]\label{alg:intersect}
{\bf Input:} an $n\times n$ integral matrix $A$ whose rows span the full
dimensional lattice $L$ in $\Q^n$, and an $m\times n$ matrix
$B$ whose rows span an $m$-dimensional $\Q$-subspace $W$ of $\Q^n$.\\
{\bf Output:} an $n\times n$ integral matrix whose rows span $L$,
and whose first $m$ rows span the lattice $W\cap L$.\\
\begin{enumerate}
\item Let $e_1,\ldots,e_n$ and $b_1,\ldots,b_m$ denote the rows of $A$ and
$B$ respectively. Write $b_i = \sum_{j=1}^n \beta_{ij}e_j$, and let
$B' =(\beta_{ij})$; if necessary multiply the rows of $B'$ by integers
in order to get integral entries.
\item Let $C$ be the output of Algorithm \ref{alg:saturation} with input $B'$.
\item Let $S$, $P$, $Q$ be the output of the Smith normal form algorithm with
input $C$.
\item return $Q^{-1}A$.
\end{enumerate}
\end{algorithm}

\begin{lemma}
Algorithm \ref{alg:intersect} is correct.
\end{lemma}

\begin{proof}
The idea is to use the given basis of $L$ as a basis of $\Q^n$. Let
$\psi: \Q^n\to \Q^n$ be the corresponding isomorphism. So if $v\in \Q^n$
then $\psi(v)$ is the vector that contains the coefficients of $v$ with
respect to the basis of $L$. So after the first step the
rows of $B'$ form a basis of $\psi(W)$. Of course $\psi(L) = \Z^n\subset \Q^n$.
\par
So the rows of $C$ form a basis of $\psi(W)\cap \psi(L)=\psi(W\cap L)$. 
Furthermore,
the Smith normal form $S$ of $C$ has diagonal entries equal to $1$. Therefore
the rows of $Q^{-1}$ form a basis of $\Z^n$ and the first $m$ rows form a 
$\Z$-basis of $\psi(W\cap L)$. Note that for $v\in \Q^n$ we have $\psi^{-1}(v)
=\psi(v)A$. Therefore the rows of $Q^{-1}A$ form a basis of 
$L$, and the first $m$ are a $\Z$-basis of $W\cap L$. 
\end{proof}

\begin{algorithm}[$L$-complements]\label{alg:complements}
{\bf Input:} an $n\times n$ integral matrix $A$ whose rows span a 
full-dimensional lattice $L$ in $V=\Q^n$; and bases of subspaces 
$V_1\subseteq V_{n-1}\subset V$.\\
{\bf Output:} an $n\times n$ integral matrix $C$ with the following properties:
\begin{itemize}
\item[-] The rows of $C$ span $L$.
\item[-] The first $s$ rows of $C$ span $L\cap V_1$ ($s=\dim V_1$).
\item[-] The first $t$ rows of $C$ span $L\cap V_{n-1}$ ($t=\dim(V_{n-1})$). 
\end{itemize}
\begin{enumerate}
\item Execute Algorithm \ref{alg:intersect} with input $A$ and a
matrix whose rows span $V_{n-1}$. Let $w_1,\ldots,w_n$ denote the
rows of the output.
\item Let $v_1,\ldots,v_s$ be the given basis of $V_1$, and write 
$v_i = \sum_{j=1}^t \alpha_{ij} w_j$. Let $A' = (\alpha_{ij})$.
\item Execute Algorithm \ref{alg:intersect} with input the $t\times t$-identity
matrix, and $A'$. Let $B$ denote the output.
\item Let $C'$ be the product of $B$ and the $t\times n$ matrix whose rows are
$w_1,\ldots,w_t$. Let $C$ be the matrix obtained from $C'$ by appending
$w_{t+1},\ldots,w_n$.
\end{enumerate}
\end{algorithm}

\begin{lemma}\label{lem:Lcomp}
Algorithm \ref{alg:complements} is correct. Let $u_1,\ldots, u_n$ denote
the rows of its output matrix $B$. Let $W_1,W_{n-1}\subset \Q^n$ be the 
subspaces spanned by $u_{s+1},\ldots,u_n$ and $u_{t+1},\ldots,u_n$ 
respectively. Then $W_1\subset W_{n-1}$ are a system of $L$-complements
to $V_1\subset V_{n-1}$.
\end{lemma}

\begin{proof}
After the first step $w_1,\ldots,w_n$ span $L$ and $w_1,\ldots,w_t$ span
$L\cap V_{n-1}$. Next we work in $V_{n-1}$ using the basis $w_1,\ldots,w_t$.
We rewrite the basis elements of $V_1$ with respect to this basis. We note
that multiplying a $v\in \Q^t$ with the matrix with rows $w_1,\ldots,w_t$ is
the reverse transformation. So 
the rows of $C'$ span $L\cap V_{n-1}$ and the first $s$ rows of $C'$ span
$L\cap V_1$. Therefore the output is correct.\par
The last statement follows directly from the definition of $L$-complements.
\end{proof}

\begin{remark}
The output of Algorithm \ref{alg:complements} has one more useful
property. Let $W_1,W_{n-1}$ be as in Lemma \ref{lem:Lcomp}.
The bases of the spaces $W_1$ and $W_{n-1}$ that are produced by this
algorithm are bases of $L\cap W_1$ and $L\cap W_{n-1}$ respectively.
\end{remark}

\begin{algorithm}[integral relations]\label{alg:intrel}
{\bf Input:} An $m\times n$-matrix $A$ with rational coefficients.\\
{\bf Output:} an  $m\times n$ integral matrix whose rows are a  
basis of the lattice
$$\Lambda = \{ (e_1,\ldots,e_m)\in \Z^m \mid \sum_{i=1}^m e_i a_i \in \Z^n\},$$
where $a_1,\ldots,a_m$ are the rows of $A$.\\
\begin{enumerate}
\item Let $M$ be the matrix obtained by appending the $n\times n$-identity
matrix at the bottom of $A$.
\item Let $v_1,\ldots,v_m$ be a basis of the space $\{ v\in \Q^{m+n} \mid
vM=0\}$. If necessary multiply each $v_i$ by an integer to ensure that
it has integral coefficients.
\item  Let $B$ be the output of the saturation algorithm (Algorithm
\ref{alg:saturation}) applied to the matrix with the $v_i$ as rows.
\item Output the rows of $B$ with the last $n$ coefficients deleted. 
\end{enumerate}
\end{algorithm}

\begin{lemma}\label{lem:intrel}
Algorithm \ref{alg:intrel} is correct.
\end{lemma}

\begin{proof}
Note that the matrix $M$ has rank $n$; therefore in Step 2. we find $m$
linearly independent basis vectors. Now set 
$$\Lambda' = \{ e=(e_1,\ldots,e_{m+n})\in \Z^{m+n} \mid eM = 0\}.$$
Then $(e_1,\ldots,e_{m+n})\mapsto (e_1,\ldots,e_m)$ is a bijection
$\Lambda\to\Lambda'$. Now after Step 3. $B$ is a basis of $\Lambda'$.
We conclude that the output is a basis of $\Lambda$. 
\end{proof}

\section{The main algorithm}\label{sect:main}

Now we return to our initial problem. Let $G\subset \GL_m(\C)$ be a
unipotent algebraic group defined over $\Q$. Set $V=\Q^m$ and let
$L$ be a full-dimensional lattice in $V$. The problem is to compute a 
finite set of generators of the group $G_L = \{ g\in G(\Q) \mid
g(L) = L\}$. \par
We assume that $G$ is defined as subset of $\GL_m(\C)$ by polynomial equations
that have coefficients in $\Q$. Then as a first step we find the Lie algebra
$\g$ of $G$ as follows. First we compute a set $\mathcal{S}$
of generators for the radical of the 
ideal generated by the polynomials that define $G$ (cf. \cite{beckwsp}).
Then we obtain the Lie algebra by differentiating the elements of $\mathcal{S}$.
\par
The second step will be to compute a flag $0=V_0 < V_1 < \cdots <V_n=V$
of $V$ for the action of $\g$. This is done by straightforward linear algebra:
$V_1$ is the space killed by all elements of $\g$; $V_2/V_1$ is the subspace
of $V/V_1$ that is killed by all elements of $\g$, and so on. \par
Now we have the input for our main algorithm which we now state. We use
the notation of Sections \ref{thresults}, \ref{sect:basis}.

\begin{algorithm}[Main algorithm]\label{alg:main}
{\bf Input:} a non-zero finite dimensional vector space $V$ over $\mathbb{Q}$, 
a full-dimensional lattice $L$ of $V$, a Lie subalgebra 
$\mathfrak{g} \subset \mathfrak{gl}(V)$ that is the Lie algebra of
a unipotent algebraic group $G$, and a flag 
\begin{displaymath}
0 = V_{0} < V_{1} < \cdots < V_{n} = V
\end{displaymath}
of $V$ with respect to the action of $\mathfrak{g}$.\\
{\bf Output:} a $T$-sequence for $G_L$.
\begin{enumerate}
\item If $n=1$ then return the empty set, else go to Step 2.
\item Compute the derived vector space $V^{\star}$, the derived lattice 
$L^{\star}$, the derived flag $0 = V^{\star}_{0} < V^{\star}_{1} < \cdots 
< V^{\star}_{n-1} = V^{\star}$.
\item Compute the kernel $\mathfrak{n}$ and the image $\mathfrak{q}$ of the 
derived action of $\mathfrak{g}$ on $V^{\star}$, together with the projection 
$\mathrm{d} \pi : \mathfrak{g} \rightarrow \mathfrak{q}$.
\item Apply the algorithm recursively to the vector space $V^{\star}$, the 
lattice $L^{\star}$, the Lie algebra $\mathfrak{q}$ and the derived flag. 
Denote by $q_{1}, \dots, q_{k}$ the result.
\item Compute a system $W_{n-1} \subseteq W_{1}$ of $L$-complements to 
$V_{1} \subseteq V_{n-1}$, (Algorithm \ref{alg:complements})
the induced lattice $\Gamma$ and the induced error map $\epsilon : 
\mathrm{End}(V) \to \mathrm{Lin}(W_{n-1}, V_{1})$.
\item Compute a basis $x_{1}, \dots, x_{l}$ of $\mathfrak{n}_{L}$ and set 
$n_i = \exp(x_i)$  for $1\leq i\leq l$.
\item For $1\leq i\leq k$ compute a preimage $x$ of $\log(q_{i})$ through 
$\mathrm{d} \pi$ and set $g_{q_{i}} = \exp(x)$.
\item Compute the image $W'$ of $\mathfrak{n}$ through $\epsilon$.
\item Compute a basis $\mathcal{W} = \{ w_{1}, \dots, w_{k} \}$ in Hermite 
normal form of the relation lattice of the elements $\epsilon(g_{q_{i}}) + 
\Gamma + W'$ in $\frac{\textrm{Lin}(W_{n-1}, V_{1})}{\Gamma + W'}$
for $i=1,\ldots,k$.
\item For each $w_{i}$ in $\mathcal{W}$ do the following
\begin{enumerate}
\item Write $w_i = (e_1^{(i)},\ldots,e_k^{(i)})$ and set 
$$g_{w_{i}}= g_{q_1}^{e_1^{(i)}}\cdots g_{q_k}^{e_k^{(i)}}.$$ 
\item \label{findw} Compute $v_{w_{i}} \in W'$ and $\gamma_{w_{i}} \in \Gamma$ such that 
$v_{w_{i}} + \gamma_{w_{i}} = \epsilon(g_{w_{i}})$.
\item Compute the preimage $n_{w_{i}}$ of $-v_{w_{i}}$ through 
$\epsilon : \mathfrak{n} \rightarrow \textrm{Lin}(W_{n-1}, V_{1})$.
\item Compute $g_{i} = g_{w_{i}} \cdot \textrm{exp}(n_{w_{i}})$.
\end{enumerate}
\item Return $g_{1}, \dots, g_{k}, n_{1}, \dots, n_{l}$.
\end{enumerate}
\end{algorithm}

We start by commenting on the computability of various steps.
\begin{enumerate}
\item[5.] We note that the output of Algorithm \ref{alg:complements}
contains bases of $L\cap W_{n-1}$ and $L\cap V_1$. We uses these bases
to represent an element of $\mathrm{Lin}(W_{n-1},V_1)$ as an 
$s\times t$-matrix (where $s=\dim(V_1)$, $t=\dim(W_{n-1})$). Then
a $\Z$-basis of $\Gamma$ is the set of elementary $s\times t$-matrices, which
have one coefficient equal to $1$, and all other coefficients equal to $0$.
Computing $\epsilon(a)$ for an $a\in \mathrm{End}(V)$ is standard linear 
algebra. Indeed, for $a\in \mathrm{End}(V)$ and $w\in W_{n-1}$ write 
$aw = v_1+w_1$, where $v_1\in V_1$ and $w_1\in W_1$. Then $\epsilon(a)w= v_1$.
\item[6.] Here we first compute a basis of the space  
$\epsilon(\mathfrak{n})\subset \mathrm{Lin}(W_{n-1},V_1)$. Using Algorithm 
\ref{alg:intersect} we find a $\Z$-basis of $\Gamma \cap 
\epsilon(\mathfrak{n})$. The inverse images of the basis elements under
$\epsilon$ are then a basis of $\mathfrak{n}_L$ (Proposition \ref{prop:ninj}).
\item[8.] In Step 6. we already obtained a basis of $W'$, and a basis
$\gamma_1,\ldots,\gamma_{st}$ of $\Gamma$ such that $\gamma_1,\ldots,\gamma_m$
are a basis of $W'\cap \Gamma$. Let $N=st-m$ and let $\psi : \mathrm{Lin}(
W_{n-1},V_1)\to \Q^N$ defined as follows. For $\gamma\in \mathrm{Lin}(W_{n-1},
V_1)$ write $\gamma =\sum_i c_i \gamma_i$; then $\psi(\gamma) = (c_{m+1},\ldots,
c_{st})$. Then $\psi$ is linear and $\gamma\in  \Gamma + W'$ if and only if
$\psi(\gamma) \in  \Z^{N}$. Set $u_i = \epsilon(g_{q_i})$ for $1\leq i\leq k$. 
We want to compute a $\Z$-basis of the lattice
$$\Lambda = \{ (e_1,\ldots,e_k)\in \Z^k \mid \sum_{i=1}^k e_iu_i \in W'+
\Gamma\}.$$
Now $(e_1,\ldots,e_k)$ lies in $\Lambda$ if and only if $\sum_i e_i \psi(u_i)
\in \Z^{N}$. So we get a basis of $\Lambda$ by applying Algorithm 
\ref{alg:intrel} with input the matrix with rows $\psi(u_i)$. Then we 
compute the Hermite 
normal form (cf. \cite{sims}, \S 8.1) of the basis of $\Lambda$ obtained.
\end{enumerate}

The other steps are straightforward. We now prove the correctness of the
algorithm.

\begin{theorem}[Correctness]\label{thm:correct}
Let $V$ be a non-zero finite dimensional vector space over $\mathbb{Q}$, 
$L$ a full-dimensional lattice of $V$, and $G$ a unipotent algebraic subgroup 
of $GL(V)$. Further, let $\mathfrak{g}\subset \mathfrak{gl}(V)$ be the Lie 
algebra of $G$. Then Algorithm \ref{alg:main}, with input $V$, $L$ and
$\mathfrak{g}$ and a flag
\begin{displaymath}
0 = V_{0} < V_{1} < \cdots < V_{n} = V
\end{displaymath}
of $V$ with respect to the action of $\mathfrak{g}$, returns a $T$-sequence 
for the subgroup $G_{L}$ of $\GL(V)$.
\end{theorem}

\begin{proof}
As a first thing we notice that, since $\mathfrak{g}$ is the Lie algebra 
of $G$, as seen in Section
\ref{sect:basis} the given flag is also a flag of $V$ with respect to the 
action $G$. Hence the hypothesis of Section \ref{thresults} is satisfied, 
and we 
can consider all the constructions described there. As seen in Section 
\ref{sect:basis}, $\mathfrak{n}\subset \mathfrak{gl}(V)$ is the Lie algebra
of $N$ and $\mathfrak{q} \subset \mathfrak{gl}(V^{\star})$ is the Lie algebra
of $Q$. Now we proceed by induction on the length $n$ of the flag. \par 
If $n = 1$, then every vector of $V$ is $G$-fixed, hence $G$ is the trivial 
subgroup of $GL(V)$ and $G_{L}$ is the trivial subgroup of $GL(V)$, hence the 
empty set is a $T$-sequence for $G_{L}$.\par 
Now suppose that $n \geq 2$. By the inductive hypothesis, $q_{1}, \dots, q_{k}$
is a $T$-sequence for the subgroup $Q_{L^{\star}}$ of $GL(V^{\star})$. By results 
in Section \ref{sect:basis}, the exponential map gives a bijection from 
$\mathfrak{n}$ to $N(\mathbb{Q})$ and, since by Proposition \ref{central} the 
Lie algebra $\mathfrak{n}$ is abelian, it is also a group morphism. Further, 
using Proposition \ref{commtriang} and Lemma \ref{epsgamma}, it is easily seen 
that the image of $\mathfrak{n}_{L}$ through the exponential map is $N_{L}$. 
Hence $n_{1}, \dots, n_{s}$ is a $T$-sequence for $N_{L}$, no matter in which 
order they are taken. As seen in Section \ref{sect:basis}, we have 
$\exp( \mathrm{d}\pi (x) ) = \pi( \exp(x) )$ for all $x\in \g$. So for 
$i = 1, \dots, k$ we get
\begin{equation} \label{eq:repr}
\pi(g_{q_{i}}) = q_{i}.
\end{equation}
In particular, $g_{q_{i}} \in G_{L^{\star}}$. By commutativity of the 
diagram in Corollary \ref{commdiag},
\begin{displaymath}
\Psi (q_{i}) = \epsilon(g_{q_{i}}) + \Gamma + W
\end{displaymath}
and, due to Proposition \ref{commtriang}, $W'$ is equal to the image $W$ of 
$N(\mathbb{Q})$ through $\epsilon$. Hence $\mathcal{W}$ is a basis in Hermite 
normal form for the relation lattice of $\Psi(q_{i}), i = 1, \dots, k$ in 
$\frac{\textrm{Lin}(W_{n-1}, V_{1})}{\Gamma + W}$. For $1\leq i\leq k$ set
$$ h_i = q_1^{e_1^{(i)}}\cdots q_k^{e_k^{(i)}},$$
where $w_i = (e_1^{(i)},\ldots,e_k^{(i)})$ is as in the algorithm.
By Lemma \ref{lemma:Eick}, the ordered set $h_{1}, \dots, h_{k}$ is a 
$T$-sequence for the kernel of $\Psi$. Due to \eqref{eq:repr}, we have that 
$g_{w_{i}}$ is an element of $G_{L^{\star}}$ satisfying
\begin{displaymath}
\pi(g_{w_{i}}) = h_{i}
\end{displaymath}
for $1\leq i\leq k$. 
Again by commutativity of the diagram in Corollary \ref{commdiag}, we have 
that $\epsilon(g_{w_{i}}) \in \Gamma + W = \Gamma + W'$ (hence Step \ref{findw} 
makes sense). Since $\textrm{exp}(n_{w_{i}}) \in N(\mathbb{Q})$, we have 
$g_{i} \in G_{L^{\star}}$ and
\begin{displaymath}
\pi(g_{i}) = h_{i}.
\end{displaymath}
Further, using Propositions \ref{action} and \ref{commtriang}, we obtain that
\begin{displaymath}
\epsilon(g_{i}) = \epsilon(g_{w_{i}}) + \epsilon(n_{w_{i}}) = v_{w_{i}} - v_{w_{i}} 
+ \gamma_{w_{i}} \in \Gamma
\end{displaymath}
hence by Lemma \ref{epsgamma} we have that $g_{i} \in G_{L}$. Thus by 
Proposition \ref{kerpsi} the ordered set $g_{1} N_{L}, \dots, g_{k} N_{L}$ is 
a $T$-sequence for $\frac{G_{L}}{N_{L}}$. Using Proposition \ref{central} 
we get that $N$ is central in $G$, hence 
$N_{L}$ is central in $G_{L}$. Therefore we finally obtain that $g_{1}, \dots, 
g_{k}, n_{1}, \dots, n_{l}$ is a $T$-sequence for $G_L$.
\end{proof}

\begin{corollary}
Let the notation be as in Theorem \ref{thm:correct}. Let $\mathcal{H}$ denote
the $T$-group generated by the output of Algorithm
\ref{alg:main}. Then the Hirsh length of $\mathcal{H}$ is equal to the dimension
of $G$. Moreover, the Lie algebra of the radicable hull of $\mathcal{H}$ is
isomorphic over $\Q$ to $\g$.
\end{corollary}

\begin{proof}
We use the notation of the proof of Theorem \ref{thm:correct}. By induction,
$q_1,\ldots,q_k$ is a $T$-sequence for $Q_{L^\star}$. This group has dimension
equal to $\dim G - \dim N = \dim G - l$. Therefore the $T$-sequence output by
the algorithm is of length equal to $\dim G$. \par
Let $h_1,\ldots,h_r$ denote the $T$-sequence output by the algorithm.
Due to \cite{segal} Chapter 6, the Lie algebra of the radicable hull
of $\mathcal{H}$ is isomorphic to the Lie algebra spanned by $\log( h_i )$.
But the latter one is $\g$. 
\end{proof}

\begin{remark}
The algorithm can be slightly modified in order to compute even a finite 
presentation of $G_{L}$. To show how it can be done, we need to introduce some 
notation first.\par 
Let $\mathcal{G}$ be a $T$-group and $g_{1}, \dots, g_{n}$ a $T$-sequence for 
$\mathcal{G}$. We 
say that a word $w$ in $g_{1}, \dots, g_{n}$ is normal if it is of the form 
$g_{1}^{e_{1}} \cdots g_{n}^{e_{n}}$ for some $e_{i} \in \mathbb{Z}$. If this is 
the case, the depth of $w$ is the minimum $i$ such that $e_{i}$ is non-zero. 
It is an easy induction to prove that every element $g$ of $\mathcal{G}$ can be written 
in a unique way as a normal word $w$ in $g_{1}, \dots, g_{n}$. Also, 
$g \in \mathcal{G}_{i} = \langle g_{i}, \dots, g_{n} \rangle$ if and only if 
$w$ has 
depth at least $i$. Since $g_{1}, \dots, g_{n}$ is a $T$-sequence, for every 
$1 \leq i < j \leq n$ we have that $[g_{i}, g_{j}] \in \mathcal{G}_{j+1}$, 
hence there 
exist a unique normal word $w_{[g_{i}, g_{j}]}$ (of depth at least $j+1$) such 
that
\begin{displaymath}
[g_{i}, g_{j}] = w_{[g_{i}, g_{j}]}.
\end{displaymath}
It is well known that
\begin{displaymath}
\left \langle g_{1}, \dots, g_{n} \ \vline \ [g_{i}, g_{j}] = w_{[g_{i}, g_{j}]} 
\ \mathrm{for} \ 1 \leq i < j \leq n \right \rangle
\end{displaymath}
is a finite presentation for $G$. We call it the standard presentation for 
$\mathcal{G}$ with respect to $g_{1}, \dots, g_{n}$.\par 
Recall that the algorithm as it stands computes a $T$-sequence 
$g_{1}, \dots, g_{k}, n_{1}, \dots, n_{l}$ for $G_{L}$. Now we want to sketch 
how it can be modified in order to compute the standard presentation for 
$G_{L}$ with respect to such a $T$-sequence. Of course it will be enough to 
show how to compute the normal words of the forms $w_{[g_{i}, g_{j}]}$, 
$w_{[g_{i}, n_{j}]}$ and $w_{[n_{i}, n_{j}]}$ (for any suitable choice of the 
indexes $i$ and $j$). The proof of Theorem \ref{thm:correct} shows that 
$N_{L}$ is central in $G_{L}$. Hence it follows at once that
\begin{displaymath}
w_{[g_{i}, n_{j}]} = w_{[n_{i}, n_{j}]} = 1.
\end{displaymath}
So the only hard part is to compute the words of the form $w_{[g_{i}, g_{j}]}$. 
To this end we can suppose by inductive hypothesis that Step $4$ of the 
algorithm, along with a $T$-sequence $q_{1}, \dots, q_{k}$ for $Q_{L^{\star}}$, 
produces also the standard presentation
\begin{displaymath}
\left \langle q_{1}, \dots, q_{k} \ \vline \ [q_{i}, q_{j}] = w_{[q_{i}, q_{j}]} 
\ \mathrm{for} \ 1 \leq i < j \leq k \right \rangle
\end{displaymath}
of $Q_{L^{\star}}$ with respect to $q_{1}, \dots, q_{k}$. Recall that, using the 
notations of the proof of Theorem \ref{thm:correct}, $h_{1}, \dots, h_{k}$ is 
a $T$-sequence for the kernel of $\Psi$. Further, the proof even shows that we 
can effectively write its elements as normal words in $q_{1}, \dots, q_{k}$ 
(of increasing depth). There are algorithms known for computing a standard
presentation of a subgroup of a $T$-group (cf. \cite{sims};
an implementation of the algorithms for this purpose is available in
the {\sf GAP}4 package ``polycyclic'', \cite{polycyclic}). So we 
can compute the standard presentation
\begin{displaymath}
\left \langle h_{1}, \dots, h_{k} \ \vline \ [h_{i}, h_{j}] = w_{[h_{i}, h_{j}]} 
\ \mathrm{for} \ 1 \leq i < j \leq k \right \rangle
\end{displaymath}
for the kernel of $\Psi$ with respect to $h_{1}, \dots, h_{k}$. Since every 
$w_{[h_{i}, h_{j}]}$ is a normal word in $h_{1}, \dots, h_{k}$, we can evaluate 
it in $g_{1}, \dots, g_{k}$ (that it to say, substituting every $h_{i}$ with 
$g_{i}$). In this way we obtain a normal word $u_{[g_{i}, g_{j}]}$ in $g_{1}, 
\dots, g_{n}$. Since $\pi : G_{L} \rightarrow K$ sends $g_{i}$ in $h_{i}$ and 
has kernel $N_{L}$, we obtain that
\begin{displaymath}
n = [g_{i}, g_{j}] u_{[g_{i}, g_{j}]}^{-1} \in N_{L}
\end{displaymath}
Since $n_{1}, \dots, n_{s}$ is a $T$-sequence for $N_{L}$, we can write $n$ as 
a normal word $v$ in $n_{1}, \dots, n_{s}$. This can be done effectively, being 
equivalent to write an element of a lattice in terms of a (ordered) basis. 
Since $u_{[g_{i}, g_{j}]} v$ is a normal word in $g_{1}, \dots, g_{k}, n_{1}, 
\dots, n_{l}$, it now follows easily that
\begin{displaymath}
w_{[g_{i}, g_{j}]} = u_{[g_{i}, g_{j}]} v.
\end{displaymath}
\end{remark}

\section{Practical performance}\label{sect:pract}
It is rather straightforward to see that the complexity of Algorithm
\ref{alg:main} is exponential in the length of the flag of $V$. Indeed,
if the flag has maximal length, or in other words, $\dim V_i =i$, then
$\dim V^\star = 2n-2$. So in the worst case the dimension of the ambient
vector space is roughly doubled in each step of the recursion. \par
It turns out that the dimension of the spaces $\mathrm{Lin}( W_{n-1}, V_1)$
increases even faster (in the worst case by about a factor of $4$ each
step of the recursion). For this reason we avoid working with the
entire space $\mathrm{Lin}( W_{n-1}, V_1)$. Instead we consider the associative
algebra with one $A\subset \mathrm{End}(V)$ generated by the elements of the 
Lie algebra $\g$. Let $U=\epsilon(A)\subset \mathrm{Lin}(W_{n-1}, V_1)$. 
Then $\epsilon(G)\subset U$, so we can work with the space $U$ instead of
$\mathrm{Lin}(W_{n-1},V_1)$. In fact we choose to work with a potentially
somewhat bigger space, namely the subspace of $\mathrm{Lin}(W_{n-1},V_1)$
consisting of the matrices that have nonzero entries only in those positions
for which there are elements in $U$ that have nonzero entries in those 
positions. This space has the advantage that its intersection with the
lattice $\Gamma$ is easily computed. Furthermore, in practice it is only
slightly bigger than $U$.\par
We have implemented the algorithm in the language of {\sf GAP}4
(\cite{gap4}). We use two series of Lie algebras in $\gl_n(\Q)$ to generate
test inputs to the algorithm. The terms of the first series are denoted
$\g_n$, which consists of the elements
\begin{align*}
x_i &= e_{1,i+1} \text{ for $1\leq i\leq n-1$ }\\
x_n &= \sum_{j=2}^{n-1}e_{j,j+1}.  
\end{align*}
(here $e_{i,j}$ is the $n\times n$-matrix with a $1$ on position $(i,j)$ and
zeros elsewhere). The only nonzero commutators are $[x_i,x_n] = x_{i+1}$
for $1\leq i\leq n-2$. So $\g_n$ is of dimension $n$ and of nilpotency
class $n-1$. \par
The second series of Lie algebras is denoted $\h_n$, which consists of
the elements
\begin{align*}
y_1 &= \sum_{i=1}^{n-1} ie_{i,i+1}\\
y_k &= \sum_{i=1}^{n-k} e_{i,i+k} \text{ for $2\leq k\leq n-1$}.
\end{align*} 
Here the only nonzero commutators are $[y_1,y_k] = -ky_{k+1}$ for
$2\leq k\leq n-2$. So $\h_n$ is of dimension $n-1$ and of nipotency
class $n-2$. In fact, as abstract Lie algebras, $\h_n\cong \g_{n-1}$.
We note that both Lie algebras have a flag of maximal length.\par
In Table \ref{table1} we list the runing times(\footnote{The computations
were done on a 2GHz processor with 1GB of memory for {\sf GAP}.}) of the 
algorithm with input the Lie algebras $\g_n$ and $\h_n$, for $n=6,7,8,9$. 
In all cases for the lattice $L$ we have taken $\Z^n\subset \Q^n$. 

\begin{table}[htb]
\begin{center}
\begin{tabular}{|r|r|r|}
\hline
$n$ & time $\g_n$ & time $\h_n$\\
\hline
6 & 0.7 & 0.4 \\
7 & 3  &  3 \\
8 & 24 & 16  \\
9 & 204 & 133\\
\hline

\end{tabular}
\caption{Time (in seconds) for the main algorithm with input $\g_n$ and
$\h_n$.}\label{table1}
\end{center}
\end{table}

From Table \ref{table1} we see that the algorithm is efficient enough to
tackle nontrivial examples. However, the running times do confirm the 
analysis above that the complextity of the algorithm is exponential. 
Also we see that the running time is essentially determined by the length
of the flag, as $\h_n\cong \g_{n-1}$, but the algorithm needs markedly longer
for $\h_n$ than for $\g_{n-1}$.

\end{document}